\documentclass[12pt]{amsart}
\usepackage{amssymb}

\newlength{\horzstretch}
\newlength{\vertstretch}
\theoremstyle{plain}
\newtheorem{lemma}{Lemma}[section]
\newtheorem{proposition}[lemma]{Proposition}
\newtheorem{theorem}[lemma]{Theorem}

\newtheorem{conjecture}[lemma]{Conjecture}

\theoremstyle{definition}

\newtheorem{example}[lemma]{Example}

\newcommand{\ggcaffil}[1]{\dedicatory{\textup{\larger{#1}}}}

\newcommand{\ggcqedsymbol}{$\square$}
\newcommand{\ggcqed}{\hbox{}\nobreak\hbox{\quad\ggcqedsymbol}}
\newcommand{\ggcendpf}{\ggcqed}
\newcommand{\ggcnopf}{\ggcqed}
\newcommand{\ggcendconj}{\ggcqed}
\newcommand{\ggcenddef}{\ggcqed}
\newcommand{\ggcendeg}{\ggcenddef}

\newcommand{\sbset}{\subseteq}
\newcommand{\mtset}{\varnothing}
\newcommand{\setmin}{-}
\newcommand{\bddots}{\mathinner{%
  \mkern1mu\raise1pt\hbox{.}%
  \mkern2mu\raise4pt\hbox{.}%
  \mkern2mu\raise7pt\hbox{.}%
  \mkern1mu}}
\newcommand{\ggcitemlimits}{\limits}
\newcommand{\ggcinlinelimits}{}

\begin{document}
\title[Row-Latin Rectangles]
  {A Matroid Generalization of a Result on Row-Latin Rectangles}
\author{Glenn G. Chappell}
\ggcaffil{Department of Mathematics, Southeast Missouri
  State University}
\address{Department of Mathematics\\
  Southeast Missouri State University\\
  Cape\break Girardeau, MO 63701\\
  USA}
\email{gchappell@semovm.semo.edu}
\date{July 8, 1998}
\begin{abstract}
  Let $A$ be an $m\times n$ matrix in which the entries of
  each row are all distinct.
  Drisko~\cite{DriA97} showed that, if $m\ge2n-1$, then $A$ has
  a \emph{transversal}:
  a set of $n$ distinct entries with no two
  in the same row or column.
  We generalize this to matrices with entries in a matroid.
  For such a matrix $A$, we show that if each row
  of $A$ forms an independent set, then we can require
  the transversal to be independent as well.
  We determine the complexity of an algorithm based on the proof
  of this result.
  Lastly, we observe that $m\ge2n-1$ appears to force the existence
  of not merely one but many transversals.
  We discuss a number of conjectures related to this observation
  (some of which involve matroids and some of which do not).
\end{abstract}
\maketitle

\section{Introduction} \label{S:intro}

We define a
\emph{partial transversal} of length $k$ in a matrix $A$
to be a set of $k$ distinct entries of $A$,
no two in the same row or column.
A \emph{transversal} is a partial transversal that meets every
column.
A \emph{Latin square} of order $n$ is an $n\times n$ matrix
in which each of the rows and columns is a permutation of
$\{1,2,\dotsc,n\}$.
The existence of partial transversals in
Latin squares has been discussed in a number of
works~\cite{BVW78,DeKe91,EHNS88,ShoP82,WooD78};
see~\cite{DeKe91} for a survey.

Stein~\cite{SteS75} and Erd{\H o}s et al.~\cite{EHNS88}
investigated partial transversals
in generalizations of Latin squares.
One such generalization is a \emph{row-Latin rectangle}:
an $m\times n$ matrix in which each row is a permutation
of $\{1,2,\dotsc,n\}$.

Dillon~\cite{DilJ96} asked for the minimum
$m$ so that every $m\times n$ row-Latin rectangle has a
transversal.
Drisko~\cite{DriA97} answered Dillon's question by showing
that $m=2n-1$ suffices.
In fact, he proved this for a more general class of matrices.
A \emph{row-Latin rectangle based on $k$} is an $m\times n$
matrix with entries in $\{1,2,\dotsc,k\}$ so that no entry
appears twice in any row.

\begin{theorem}[Drisko 1997 {\cite[Thm.~1]{DriA97}}] \label{T:driskos}
Let $A$ be an $m\times n$ row-Latin rectangle based on $k$.
If $m\ge 2n-1$,
then $A$ has a transversal.\ggcnopf\end{theorem}

The following well known example,
based on~\cite[Example~1]{DriA97},
shows that the bound on $m$ in Theorem~\ref{T:driskos} is sharp.

\begin{example} \label{G:wellknown}
Let $m$, $n$ be positive integers with $m\ge n-1$.
We define $R_{m,n}$
to be an $m\times n$ matrix whose first $m-(n-1)$
rows consist of the symbols $1,2,\dotsc,n$ in order
and whose remaining $n-1$ rows
have the same symbols in the order
$2,3,\dotsc,n,1$.
Figure~\ref{F:r64} shows $R_{2,2}$, $R_{4,3}$, and $R_{6,4}$.

\begin{figure}[htbp]
\[
R_{2,2}=\begin{pmatrix}
  1&2\\
  2&1\end{pmatrix};
\qquad
R_{4,3}=\begin{pmatrix}
  1&2&3\\1&2&3\\
  2&3&1\\2&3&1\end{pmatrix};
\qquad
R_{6,4}=
\begin{pmatrix}
  1&2&3&4\\1&2&3&4\\1&2&3&4\\
  2&3&4&1\\2&3&4&1\\2&3&4&1\end{pmatrix}.
\]
\caption{The matrices $R_{2,2}$, $R_{4,3}$, and $R_{6,4}$
  of Example~\ref{G:wellknown}.
  None of these matrices has a transversal.}
\label{F:r64}
\end{figure}

The matrix $R_{2n-2,n}$ is a row-Latin rectangle
with no transversals.
To see this
assume that $R_{2n-2,n}$ has a transversal.
Without loss of generality, we may assume that this
transversal has a $1$ in the first column.
Then it cannot have a $1$ in column $n$,
so it must have an $n$ in column $n$.
Similarly, it must have an $n-1$ in column $n-1$,
an $n-2$ in column $n-2$, etc., and a $3$ in column $3$.
But this leaves no possible value in column $2$,
since we have already used each of the first $n-1$ rows;
thus, there is no transversal.\ggcendeg\end{example}

Our main result
is a matroid generalization of Theorem~\ref{T:driskos}.
Let $A$ be an $m\times n$ matrix with entries in a matroid $M$.
We define an \emph{independent partial transversal} (IPT) of
length $k$ in $A$ to be a partial transversal of length $k$
whose elements form an independent set.
An \emph{independent transversal} (IT) is an IPT that meets
every column.
Suppose that the entries of each row of $A$
are all distinct and form an independent set.
We show that if $m\ge2n-1$, then $A$ has an IT.
Theorem~\ref{T:driskos} follows
by letting $M$ be a free matroid.

In Section~\ref{S:main} we prove our main result.
Our proof can be written as an algorithm to find an IT;
in Section~\ref{S:alg} we determine the complexity of this algorithm.
In Section~\ref{S:open} we observe that $m\ge 2n-1$ appears to
force the existence of not merely one but many transversals.
We discuss a number of conjectures and examples
stemming from this observation.

\section{The Main Result} \label{S:main}

\begin{theorem} \label{T:main}
Let $A$ be an $m\times n$ matrix with entries in a matroid $M$.
Suppose that the set of entries of each row of $A$
forms an independent set of size $n$ in $M$.
If $m\ge 2n-1$,
then $A$ has an IT.
\end{theorem}

\begin{proof}
Our proof is based on a simplification of Drisko's proof of
Theorem~\ref{T:driskos}~\cite[Thm.~1]{DriA97}.

Let $m$, $n$, $A$, and $M$ be as in the statement
of the theorem.
For a set $S$ in $M$, $r(S)$ denotes the rank of $S$,
and $\sigma(S)$ denotes the span of $S$.
We assume that the entries of $A$ are all distinct;
if not, we can blow up parallel classes.
Given a set $S$ of entries of $A$ and an entry $a$, we write
$S+a$ for $S\cup\{a\}$ and $S-a$ for $S\setmin\{a\}$.

It suffices to prove the result when $m=2n-1$.
The $n=1$ case is trivial; we proceed by induction on $n$.
We will prove the $n=2$ case as we set up the
induction.

We first name a number of entries of $A=(a_{i,j})$.
We let $b_1$ and $b_2$ denote $a_{1,1}$ and $a_{2,1}$, respectively.
For $1\le i\le n-1$, $c_i$ denotes $a_{i+2,i+1}$.
For $1\le i\le n-3$, $d_i$ denotes $a_{n+i+1,n-i+1}$.
Note that if $n=2,3$, then we do not define any $d_i$'s.
See Figure~\ref{F:bigmat}.

\begin{figure}[htbp]
\newcommand{\tallstrut}{\vbox{\mathstrut\kern1ex}}
\newcommand{\evs}[1]{%
  \setbox0=\hbox{$c_1$}%
  \hskip5pt\hbox to\wd0{\mathstrut\hss$#1$\hss}\hskip5pt%
  }
\newcommand{\levs}[1]{%
  \setbox0=\hbox{$c_1$}%
  \hbox to\wd0{\mathstrut\hss$#1$\hss}\hskip5pt%
  }
\newcommand{\revs}[1]{%
  \setbox0=\hbox{$c_1$}%
  \hskip5pt\hbox to\wd0{\mathstrut\hss$#1$\hss}%
  }
\newcommand{\xevs}{%
  \setbox0=\hbox{$c_1$}%
  \setbox2=\hbox{$c_{n-1}$}%
  \hbox{\kern0.5\wd2\kern-0.5\wd0}%
  }
\[
\arraycolsep=0pt
\overset{
  \displaystyle
  \begin{pmatrix}
  \levs{b_1}& \\
  \levs{b_2}& \\
  &\revs{c_1}
  \end{pmatrix}
  }
{n=2\tallstrut}
\qquad
\overset{
  \displaystyle
  \begin{pmatrix}
  \levs{b_1}&& \\
  \levs{b_2}&& \\
  &\evs{c_1}& \\
  &&\revs{c_2} \\
  \levs{\mathrm{X}}&\evs{\mathrm{X}}&\revs{\mathrm{X}}
  \end{pmatrix}
  }
{n=3\tallstrut}
\qquad
\overset{
  \displaystyle
  \begin{pmatrix}
  \levs{b_1}&&&&&& \\
  \levs{b_2}&&&&&& \\
  &\evs{c_1}&&&&& \\
  &&\evs{c_2}&&&& \\
  &&&\evs{c_3}&&& \\
  &&&&\evs{\ddots}&& \\
  &&&&&\evs{c_{n-2}}& \\
  &&&&&&\revs{c_{n-1}}&\xevs \\
  &&&&&&\revs{d_1} \\
  &&&&&\evs{d_2}& \\
  &&&&\evs{\bddots}&& \\
  &&&\evs{d_{n-3}}&&& \\
  \levs{\mathrm{X}}&\evs{\mathrm{X}}&\evs{\mathrm{X}}&&&&
  \end{pmatrix}
  }
{n\ge4\tallstrut}
\]
\caption{Various named entries for the proof of Theorem~\ref{T:main};
  X's are in positions $(2n-1,1)$, $(2n-1,2)$, and $(2n-1,3)$.}
\label{F:bigmat}
\end{figure}

If, for each column of $A$, the entries of the column are all
parallel, then the main diagonal of $A$ is an IT.
Thus,
we may assume that some column of $A$ contains two nonparallel entries.
Since we can permute the rows and columns of $A$ to place these
two entries in positions $(1,1)$ and $(2,1)$, we may assume
that the two nonparallel entries are $b_1$ and $b_2$.

Throughout this proof, we will make assumptions similar to
that above based on the fact that we could permute rows
and columns to put certain values in the required positions.
Some of these permutations may move the entries of
previously defined sets;
however, we will always ensure that the set of positions
occupied by the elements of each such set does not change.
For example, if $S=\{c_1,c_2\}$, then the permutation that
transposes rows 3 and 4 and transposes columns 2 and 3 moves
elements of $S$.
However, the \emph{set} of positions occupied by the elements of $S$
does not change.

Now, $\{b_1,b_2\}$ is an independent set of size $2$,
and $\{c_1\}$ is an independent set of size $1$.
We may augment
$\{c_1\}$ from $\{b_1,b_2\}$ to produce an IPT of length 2;
this is an  IT if $n=2$.
Thus, we may assume $n\ge3$.

If we delete the first two rows and the first column of $A$, we obtain a
$(2n-3)\times(n-1)$ matrix with entries in $M$
in which the elements of each row form an independent set.
By the
induction hypothesis, this matrix has an IT $P_1$,
which is an IPT of length $n-1$ in $A$.
The set $P_1$
meets neither the first two rows nor the first column of $A$.
Permuting rows
and columns, we may assume that $P_1=\{c_1,c_2,\dotsc,c_{n-1}\}$.

For each set $S$ of entries of $A$ with $c_1,c_2\in S$ and
$b_1,b_2\not\in S$, we define $S':=S-c_1+b_1$
and $S'':=S-c_2+b_2$.

If $P_1+b_1$ is independent in $M$, then $P_1+b_1$ is an
IT, and we are done.
A similar argument applies to $P_1+b_2$, and so
we may assume that $P_1+b_1$ and $P_1+b_2$ are dependent.
Thus,
there is a unique circuit $C_1$ with $b_1\in C_1\sbset P_1+b_1$.
Similarly,
there is a unique circuit $C_2$ with $b_2\in C_2\sbset P_1+b_2$.
Consider $(C_1\cup C_2)\cap P_1$.
If this set contains only one element, then,
by circuit elimination, $(C_1\cup C_2)\setmin P_1$ contains
a circuit.
However, $(C_1\cup C_2)\setmin P_1=\{b_1,b_2\}$ is independent.
Thus, there exist $c_i,c_j\in P_1$ with $i\ne j$,
$c_i\in C_1$, and $c_j\in C_2$.
Permuting rows and columns in such a way as not to change the
set of positions occupied by $P_1$,
we may assume that $i=1$ and $j=2$.
Now, $P_1+b_1$ is a dependent set containing a unique circuit
$C_1$, which contains $c_1$.
Thus, $P_1+b_1-c_1=P'_1$ is an IPT.
Similarly, $P_1+b_2-c_2=P''_1$ is an IPT.

Since $P_1+b_1$ is a dependent set of
rank $n-1$ containing
two independent sets $P_1$ and $P'_1$, both of size $n-1$,
we must have $\sigma(P_1)=\sigma(P'_1)$.
Similarly, $\sigma(P_1)=\sigma(P''_1)$.

Now we have constructed $P_1$ and determined some
of its properties.
Based on the assumption that $A$ does not have an IT,
we show that there is a permutation of the rows and columns of $A$
for which the following claim holds.
We will then use the claim to verify that $A$ does have an
IT, thus proving the theorem.

\medskip\noindent\emph{Claim.}
For $1\le k\le n-2$,
there exists a set $P_k$
of entries of $A$ such that
\begin{enumerate}
\item $\{c_1,\dotsc,c_{n-k}\}\sbset P_k\sbset
\{c_1,\dotsc,c_{n-1},d_1,\dotsc,d_{k-1}\}$
  \label{I:psbset}
(where $\{d_1,\dotsc,d_{k-1}\}=\mtset$ if $k=1$),
\item $P_k$, $P'_k$, and $P''_k$ are IPTs of length $n-1$ in $A$,
  \label{I:piptr}
and
\item $\sigma(P_k)=\sigma(P'_k)=\sigma(P''_k)$.
  \label{I:psigma}
\end{enumerate}
\newcounter{ggcsaveenumi}
\setcounter{ggcsaveenumi}{\value{enumi}}
Furthermore,
\begin{enumerate}
\setcounter{enumi}{\value{ggcsaveenumi}}
\item $r\left[\bigcap\ggcitemlimits_{i=1}^k\sigma(P_i)\right]=n-k$.
  \label{I:pinter}
\end{enumerate}

\medskip\noindent\emph{Proof of Claim.}
We have already defined $P_1$.
By the earlier discussion, the claim holds for
$k=1$.
If $n=3$, then the claim is proven;
we may assume $n\ge4$.
We proceed by induction on $k$.
Let $1\le k\le n-3$, and suppose the claim
holds for $1,\dotsc,k$.
We wish to define $P_{k+1}$ so that the claim
holds for $k+1$.

Since $r\left[\bigcap\ggcinlinelimits_{i=1}^k\sigma(P_i)\right]=n-k$,
there is an entry $x$ of row $n+k+1$
that lies in the first $n-k+1$ columns
and is not in
$\bigcap\ggcinlinelimits_{i=1}^k\sigma(P_i)$.
There must exist $t$, $1\le t\le k$, with
$x\not\in\sigma(P_t)$.
Now, $P_t$, $P'_t$, and $P''_t$ have the same span
and miss columns $1$, $2$, and $3$, respectively.
Thus, if $x$ lies in one of the first 3 columns, then
one of $P_t+x$, $P'_t+x$, or $P''_t+x$ is an IT,
and so we may assume that
$x$ does not lie in the first three columns.
Permuting rows and columns, we may assume that
$x=d_k$;
note that we can choose the permutations in
such a way that
the sets of positions occupied by
$P_i$, $P'_i$, and $P''_i$ remain unchanged,
for $1\le i\le k$.

Let $P_{k+1}=P_t-c_{n-k}+d_k$.
Then $P_{k+1}$, $P'_{k+1}$, $P''_{k+1}$ are IPTs,
since $d_k\not\in\sigma(P_t)=\sigma(P'_t)=\sigma(P''_t)$.
We see that Statements~\ref{I:psbset} and \ref{I:piptr} in the
claim hold for $k+1$.
It remains to verify Statements~\ref{I:psigma} and \ref{I:pinter}.

We may assume that $P_{k+1}+b_1$ is dependent;
otherwise, it is an IT.
So $P_{k+1}+b_1$ is a dependent set of rank $n-1$
containing two independent sets of size $n-1$:
$P_{k+1}$ and $P'_{k+1}$.
We conclude that $\sigma(P_{k+1})=\sigma(P'_{k+1})$.
Similarly, $\sigma(P_{k+1})=\sigma(P''_{k+1})$,
and so Statement~\ref{I:psigma} holds.

To see Statement~\ref{I:pinter}, note that
$\{c_1,c_2,\dotsc,c_{n-k}\}$ forms an independent set of
rank $n-k$.
Hence,
\begin{alignat*}{2}
\bigcap_{i=1}^{k+1}\sigma(P_i)
  &=\left[\bigcap_{i=1}^k\sigma(P_i)\right]\cap\sigma(P_{k+1})
  &&\\
&=\sigma(c_1,c_2,\dotsc,c_{n-k})\cap\sigma(P_{k+1})
  &&\quad\text{by~\ref{I:psbset} and~\ref{I:pinter}}\\
&=\sigma(c_1,c_2,\dotsc,c_{n-k-1})
  &&\quad\text{since $d_k\not\in\sigma(c_1,\dotsc,c_{n-k})$.}
\end{alignat*}
Thus, Statement~\ref{I:pinter} holds.
By induction, the claim is proven.
\medskip

Now we use the claim to prove the theorem.
We apply the claim with $k=n-2$.
Since
\[
r\left[\bigcap_{i=1}^k\sigma(P_i)\right]
  =n-k=2,
\]
there is an entry $x$ of row $2n-1$
that lies in the first $3$ columns
and is not in
$\bigcap\ggcinlinelimits_{i=1}^k\sigma(P_i)$.
There must exist $t$, $1\le t\le k$, with
$x\not\in\sigma(P_t)$.
Since $x$ lies in one of the first 3 columns,
one of $P_t+x$, $P'_t+x$, or $P''_t+x$ is an IT,
and the theorem is proven.\ggcendpf\end{proof}

\section{Algorithmic Complexity} \label{S:alg}

Drisko~\cite{DriA97} noted that his proof of
Theorem~\ref{T:driskos}~\cite[Thm.~1]{DriA97} would lead to a
recursive algorithm to find a transversal.
He calculates the time complexity of this algorithm to be
\[
\frac{1}{6}n^4+n^3-\frac{13}{6}n^2+3n-2=O(n^4),
\]
where $n$ is the number of columns of the given matrix.

Our proof of Theorem~\ref{T:main} is based on the above-mentioned
proof of Drisko;
like that proof, it can be phrased as a recursive algorithm.
We briefly examine the complexity of this algorithm.

The algorithm follows the steps of the proof,
stopping if it finds an IT.
It is given an $m\times n$ matrix ($m\ge 2n-1$)
with entries in a matroid $M$, in which the entries of
each row form an independent set of size $n$ in $M$.
It returns an IT.
We assume
that $m=2n-1$,
and
that the matroid $M$ is accessed via an
independence oracle.

The algorithm begins by checking whether the entries of each
column are all parallel.
Thus, for each of $n$ columns, we search the entries in rows
$2$ through $2n-1$ for an entry that is not parallel to the
entry in row $1$.
This requires $O(n^2)$ calls to the
independence oracle.

We permute rows and columns so that two nonparallel
entries are in the positions of $b_1$ and $b_2$ and
recursively call the algorithm to find the IPT
$P_1$.
The row and column permutations are not a major factor
in the complexity of the algorithm;
the recursive call will be discussed later.

Next we find which entries should be $c_1$ and $c_2$.
This requires $O(n)$ calls to the oracle.

Lastly, we construct the IPTs
$P_2$ through $P_{n-2}$ and return the IT.
For each of these IPTs,
we find the entry that will be $d_i$ using $O(n)$ calls to
the oracle;
then we determine which of the previous IPTs
to make the new transversal out of using
another $O(n)$ calls to the oracle.
Thus, in constructing these IPTs, we make $O(n^2)$ calls
to the oracle.

Overall, running the algorithm requires $O(n^2)$ calls to the
oracle and $1$ recursive call.
Since the recursion has depth $O(n)$,
the entire algorithm makes $O(n^3)$ calls to the oracle.

In the case when $M$ is a free matroid, a multiset is independent
precisely when its elements are all distinct;
this can be determined for a multiset of $n$ elements in
$O(n)$ time.
Thus, the algorithm corresponding to Drisko's result
(Theorem~\ref{T:driskos}) runs in $O(n^4)$ time;
this agrees with Drisko's calculation.

\section{Open Problems} \label{S:open}

Recall the matrices $R_{m,n}$ of Example~\ref{G:wellknown}.
We noted that
$R_{2n-2,n}$ has no transversals.
By Theorem~\ref{T:driskos},
$R_{2n-1,n}$ has a transversal.
For example, we can see in Figure~\ref{F:rx3} that the main diagonal
of $R_{5,3}$ is a transversal.

\begin{figure}[htbp]
\[
R_{4,3}=\begin{pmatrix}
  1&2&3\\1&2&3\\
  2&3&1\\2&3&1\end{pmatrix};
\qquad
R_{5,3}=\begin{pmatrix}
  1&2&3\\1&2&3\\1&2&3\\
  2&3&1\\2&3&1\end{pmatrix};
\qquad
R_{6,3}=\begin{pmatrix}
  1&2&3\\1&2&3\\1&2&3\\1&2&3\\
  2&3&1\\2&3&1\end{pmatrix}.
\]
\caption{Matrices from Example~\ref{G:wellknown}:
  $R_{4,3}$, which has no transversals,
  and $R_{5,3}$ and $R_{6,3}$, which have many transversals.}
\label{F:rx3}
\end{figure}

However, although $R_{4,3}$ has no transversals at all,
$R_{5,3}$ has many transversals:
there are three pairwise disjoint transversals in the first three rows.
Similarly, $R_{6,3}$ has four pairwise disjoint transversals.
Generally, while $R_{m,n}$ has no transversals if $m<2n-1$,
it has $m-(n-1)$ pairwise disjoint transversals if $m\ge2n-1$.
We conjecture that this holds for more general matrices.

\begin{conjecture} \label{J:djtr}
Let $A$ be an $m\times n$ row-Latin rectangle based on $k$.
If $m\ge2n-1$, then $A$ has $m-(n-1)$ pairwise disjoint
transversals.\ggcendconj\end{conjecture}

Example~\ref{G:wellknown} shows that
Conjecture~\ref{J:djtr} is best-possible.

The matrix $R_{5,3}$ not only has three pairwise disjoint
transversals;
it has three rows that together are the union of three
transversals.
Similarly, for $m\ge2n-1$, $R_{m,n}$ has $n$ rows
that are the union of $n$ transversals.
However, more general matrices may not have this property.
Below, we construct $(n^2-1)\times n$ matrices in which the
entries of each row are all distinct,
but no $n$ rows are the union of $n$ transversals.
We conjecture that these are the largest such matrices,
that is, that $n^2$ rows force the existence of $n$
transversals whose union is $n$ rows.

\begin{example} \label{G:n2-1byn}
For each integer $n\ge2$ we define $T_n$ to
be an  $(n^2-1)\times n$ matrix as follows.
Begin with an $(n^2-1)\times(n+1)$ matrix whose rows
are $n-1$ copies of each of the cyclic permutations of
$1,2,\dotsc,n+1$.
Delete the last column of this matrix to obtain $T_n$.

Thus, the first $n-1$ rows of $T_n$ are all $1,2,3,\dotsc,n$.
The next $n-1$ rows are all $2,3,\dotsc,n,n+1$, and so on.
Each row omits exactly one element of
$\{1,2,\dotsc,n+1\}$.

The first three matrices $T_n$ are shown
in Figure~\ref{F:tns}.\ggcendeg\end{example}

\begin{figure}[htbp]
\[
T_2=\begin{pmatrix}1&2\\2&3\\3&1\end{pmatrix};
\qquad
T_3=\begin{pmatrix}
  1&2&3\\1&2&3\\
  2&3&4\\2&3&4\\
  3&4&1\\3&4&1\\
  4&1&2\\4&1&2\end{pmatrix};
\qquad
T_4=\begin{pmatrix}
  1&2&3&4\\1&2&3&4\\1&2&3&4\\
  2&3&4&5\\2&3&4&5\\2&3&4&5\\
  3&4&5&1\\3&4&5&1\\3&4&5&1\\
  4&5&1&2\\4&5&1&2\\4&5&1&2\\
  5&1&2&3\\5&1&2&3\\5&1&2&3\\\end{pmatrix}.
\]
\caption{The matrices $T_2$, $T_3$, and $T_4$ of Example~\ref{G:n2-1byn}.}
\label{F:tns}
\end{figure}

\begin{proposition}\label{P:no-nrowtr}
There do not exist $n$
rows of $T_n$ that together are the union of $n$
transversals.
\end{proposition}

\begin{proof}
Assume for a contradiction that there are $n$ rows
of $T_n$ that are the union of $n$ pairwise disjoint transversals.
Form a
matrix $A=(a_{i,j})$ with these $n$ rows, and consider the
above-mentioned transversals as transversals of $A$.
Each row of $A$ and each transversal of $A$ omit
exactly one element of $\{1,2,\dotsc,n+1\}$.
The multiset union of the
$n$ rows and the multiset union of the $n$ transversals are the same.
Equivalently, the multiset of elements omitted from the $n$
rows is equal to
the multiset of elements omitted from the $n$ transversals.
Thus, there is a
one-to-one correspondence between the set of rows and the set of
transversals so that a corresponding row-transversal pair each
omit the same number.
Number the transversals from $1$ to $n$ so that for each $i$,
transversal $i$ omits the same number as row $i$.

For the remainder of this proof, arithmetic will be modulo $n+1$.

Let $b_i$ denote the number omitted from row $i$ and transversal $i$
($1\le i\le n$), so that $a_{i,j}=b_i+j$.
Let $t_{i,j}$ denote the element of
transversal $i$ that lies in column $j$, and let $r_{i,j}$ denote the row
that this entry lies in, so that $t_{i,j}=a_{r_{i,j},j}=b_{r_{i,j}}+j$.
For any given $i$,
since row $i$ and transversal $i$ omit the same number,
we have
\[
\sum_{j=1}^na_{i,j}=\sum_{j=1}^nt_{i,j}.
\]
As stated above, $a_{ij}=b_i+j$, and $t_{ij}=b_{r_{ij}}+j$.
Subtracting,
we obtain $a_{ij}-t_{ij}=b_i-b_{r_{ij}}$, and so
\[
0=\sum_{j=1}^n(a_{i,j}-t_{i,j})=\sum_{j=1}^n(b_i-b_{r_{i,j}})
  =nb_i-\sum_{k=1}^n b_k=-b_i-\sum_{k=1}^n b_k,
\]
since $n$ is congruent to $-1$ modulo $n+1$.
Hence, we have
\[
b_i=-\sum_{k=1}^n b_k,\quad 1\le i\le n.
\]
Since the right-hand side does not depend on $i$, all of the $b_i$'s are
equal, and so all $n$ rows must omit the same number.
But by the definition of $T_n$,
at most $n-1$ rows all omit the same number.
By contradiction,
the proposition is proven.\ggcendpf\end{proof}

\begin{conjecture} \label{J:nrowtr}
Let $A$ be an $m\times n$ row-Latin rectangle based on $k$.
If $m\ge n^2$, then there exist $n$ rows of $A$ that together
are the union of $n$ transversals.
\ggcendconj\end{conjecture}

The bound on $m$
in Conjecture~\ref{J:nrowtr} is best-possible,
by Proposition~\ref{P:no-nrowtr}.

If Conjecture~\ref{J:nrowtr} is proven, then we need only verify
Conjecture~\ref{J:djtr} for $m<n^2$;
we can then prove Conjecture~\ref{J:djtr} by an inductive argument.
Hence, Conjecture~\ref{J:djtr} can be
verified for any given value of $n$ by a bounded search.

We used the matrices $R_{2n-2,n}$ of Example~\ref{G:wellknown}
to show that a $(2n-2)\times n$ matrix in which the entries of
each row are all distinct need not have a transversal.
Drisko~\cite{DriA97} conjectured that these are essentially
the only such matrices without transversals.

\begin{conjecture}[Drisko 1997 {\cite[Conj.~2]{DriA97}}] \label{J:driskos}
Let $n\ge2$.
Let $A$ be a $(2n-2)\times n$ row-Latin rectangle based on $k$.
Then either $A$ has a transversal,
or $A$ can be transformed into $R_{2n-2,n}$
by permuting rows, columns, and symbols.
\ggcendconj\end{conjecture}

The corresponding statement
for Conjectures~\ref{J:djtr} and~\ref{J:nrowtr} is false;
that is, the matrices $R_{2n-2,n}$ and $T_n$ are not unique.
For example, the left-hand matrix in Figure~\ref{F:nott3} is
a $4\times3$ row-Latin rectangle that does not have 3 pairwise
disjoint transversals.
This matrix cannot be transformed into $R_{4,3}$ by permuting
rows, columns, and symbols.
The right-hand matrix is
an $8\times3$ row-Latin rectangle in which no $3$ rows are
the union of $3$ transversals.
This matrix cannot be transformed into $T_3$ by permuting
rows, columns, and symbols.

\begin{figure}[htbp]
\[
\begin{pmatrix}
  1&2&3\\1&2&3\\
  2&1&3\\3&2&1\end{pmatrix}
\qquad\qquad
\begin{pmatrix}
  1&2&3\\1&2&3\\
  2&3&1\\2&3&1\\
  1&3&2\\1&3&2\\
  2&1&3\\2&1&3\end{pmatrix}
\]
\caption{A $4\times3$ matrix, different from $R_{4,3}$, that
  does not have 3 pairwise disjoint transversals,
  and an $8\times3$ matrix, different from $T_3$,
  in which no 3 rows are the union of 3 transversals.}
\label{F:nott3}
\end{figure}

We generalize Conjectures~\ref{J:djtr} and~\ref{J:nrowtr}
to matrices with entries in a matroid.

\begin{conjecture} \label{J:djitr}
Let $A$ be an $m\times n$ matrix with entries in a matroid $M$.
Suppose that the set of entries of each row of $A$
forms an independent set of size $n$ in $M$.
If $m\ge 2n-1$, then $A$ has $m-(n-1)$ pairwise disjoint
ITs.
\ggcendconj\end{conjecture}

\begin{conjecture} \label{J:nrowitr}
Let $A$ be an $m\times n$ matrix with entries in a matroid $M$.
Suppose that the set of entries of each row of $A$
forms an independent set of size $n$ in $M$.
If $m\ge n^2$, then there exist $n$ rows of $A$ that together
are the union of $n$ ITs.
\ggcendconj\end{conjecture}

Conjectures~\ref{J:djitr} and~\ref{J:nrowitr} imply
Conjectures~\ref{J:djtr} and~\ref{J:nrowtr}, respectively,
by letting $M$ be a free matroid.
As above,
Conjecture~\ref{J:djitr} and the bound on $m$ in
Conjecture~\ref{J:nrowitr} are best-possible.
If Conjecture~\ref{J:nrowitr} is proven, then we need only verify
Conjecture~\ref{J:djitr} for $m<n^2$.

\end{document}